\documentclass{article}
\usepackage{amssymb,amsmath,theorem,euscript}

\input{epsf}

\newcounter{sec}

\def\sm{\smallskip}

\newcounter{punct}[sec]

\def\punct{\refstepcounter{punct}{\arabic{sec}.\arabic{punct}.  }} 

\def\COUNTERS{\addtocounter{sec}{1}
              \setcounter{punct}{0}
          \setcounter{equation}{0}
          \setcounter{theorem}{0}
                  }

\newtheorem{theorem}{Theorem}[sec]
\newtheorem{proposition}[theorem]{Proposition}

\begin{document}

 \def\ov{\overline}
\def\wt{\widetilde}
 \newcommand{\rk}{\mathop {\mathrm {rk}}\nolimits}
\newcommand{\Aut}{\mathop {\mathrm {Aut}}\nolimits}
\newcommand{\Out}{\mathop {\mathrm {Out}}\nolimits}
 \newcommand{\tr}{\mathop {\mathrm {tr}}\nolimits}
  \newcommand{\diag}{\mathop {\mathrm {diag}}\nolimits}
  \newcommand{\supp}{\mathop {\mathrm {supp}}\nolimits}
  \newcommand{\indef}{\mathop {\mathrm {indef}}\nolimits}
  \newcommand{\dom}{\mathop {\mathrm {dom}}\nolimits}
  \newcommand{\im}{\mathop {\mathrm {im}}\nolimits}
 
\renewcommand{\Re}{\mathop {\mathrm {Re}}\nolimits}

\def\Br{\mathrm {Br}}

\def\SL{\mathrm {SL}}
\def\Diag{\mathrm {Diag}}
\def\SU{\mathrm {SU}}
\def\GL{\mathrm {GL}}
\def\U{\mathrm U}
\def\OO{\mathrm O}
 \def\Sp{\mathrm {Sp}}
 \def\SO{\mathrm {SO}}
\def\SOS{\mathrm {SO}^*}
 \def\Diff{\mathrm{Diff}}
 \def\Vect{\mathfrak{Vect}}
\def\PGL{\mathrm {PGL}}
\def\PU{\mathrm {PU}}
\def\PSL{\mathrm {PSL}}
\def\Symp{\mathrm{Symp}}
\def\End{\mathrm{End}}
\def\Mor{\mathrm{Mor}}
\def\Aut{\mathrm{Aut}}
 \def\PB{\mathrm{PB}}
 \def\cA{\mathcal A}
\def\cB{\mathcal B}
\def\cC{\mathcal C}
\def\cD{\mathcal D}
\def\cE{\mathcal E}
\def\cF{\mathcal F}
\def\cG{\mathcal G}
\def\cH{\mathcal H}
\def\cJ{\mathcal J}
\def\cI{\mathcal I}
\def\cK{\mathcal K}
 \def\cL{\mathcal L}
\def\cM{\mathcal M}
\def\cN{\mathcal N}
 \def\cO{\mathcal O}
\def\cP{\mathcal P}
\def\cQ{\mathcal Q}
\def\cR{\mathcal R}
\def\cS{\mathcal S}
\def\cT{\mathcal T}
\def\cU{\mathcal U}
\def\cV{\mathcal V}
 \def\cW{\mathcal W}
\def\cX{\mathcal X}
 \def\cY{\mathcal Y}
 \def\cZ{\mathcal Z}
\def\0{{\ov 0}}
 \def\1{{\ov 1}}
 \def\frA{\mathfrak A}
 \def\frB{\mathfrak B}
\def\frC{\mathfrak C}
\def\frD{\mathfrak D}
\def\frE{\mathfrak E}
\def\frF{\mathfrak F}
\def\frG{\mathfrak G}
\def\frH{\mathfrak H}
\def\frI{\mathfrak I}
 \def\frJ{\mathfrak J}
 \def\frK{\mathfrak K}
 \def\frL{\mathfrak L}
\def\frM{\mathfrak M}
 \def\frN{\mathfrak N} \def\frO{\mathfrak O} \def\frP{\mathfrak P} \def\frQ{\mathfrak Q} \def\frR{\mathfrak R}
 \def\frS{\mathfrak S} \def\frT{\mathfrak T} \def\frU{\mathfrak U} \def\frV{\mathfrak V} \def\frW{\mathfrak W}
 \def\frX{\mathfrak X} \def\frY{\mathfrak Y} \def\frZ{\mathfrak Z} \def\fra{\mathfrak a} \def\frb{\mathfrak b}
 \def\frc{\mathfrak c} \def\frd{\mathfrak d} \def\fre{\mathfrak e} \def\frf{\mathfrak f} \def\frg{\mathfrak g}
 \def\frh{\mathfrak h} \def\fri{\mathfrak i} \def\frj{\mathfrak j} \def\frk{\mathfrak k} \def\frl{\mathfrak l}
 \def\frm{\mathfrak m} \def\frn{\mathfrak n} \def\fro{\mathfrak o} \def\frp{\mathfrak p} \def\frq{\mathfrak q}
 \def\frr{\mathfrak r} \def\frs{\mathfrak s} \def\frt{\mathfrak t} \def\fru{\mathfrak u} \def\frv{\mathfrak v}
 \def\frw{\mathfrak w} \def\frx{\mathfrak x} \def\fry{\mathfrak y} \def\frz{\mathfrak z} \def\frsp{\mathfrak{sp}}
 \def\bfa{\mathbf a} \def\bfb{\mathbf b} \def\bfc{\mathbf c} \def\bfd{\mathbf d} \def\bfe{\mathbf e} \def\bff{\mathbf f}
 \def\bfg{\mathbf g} \def\bfh{\mathbf h} \def\bfi{\mathbf i} \def\bfj{\mathbf j} \def\bfk{\mathbf k} \def\bfl{\mathbf l}
 \def\bfm{\mathbf m} \def\bfn{\mathbf n} \def\bfo{\mathbf o} \def\bfp{\mathbf p} \def\bfq{\mathbf q} \def\bfr{\mathbf r}
 \def\bfs{\mathbf s} \def\bft{\mathbf t} \def\bfu{\mathbf u} \def\bfv{\mathbf v} \def\bfw{\mathbf w} \def\bfx{\mathbf x}
 \def\bfy{\mathbf y} \def\bfz{\mathbf z} \def\bfA{\mathbf A} \def\bfB{\mathbf B} \def\bfC{\mathbf C} \def\bfD{\mathbf D}
 \def\bfE{\mathbf E} \def\bfF{\mathbf F} \def\bfG{\mathbf G} \def\bfH{\mathbf H} \def\bfI{\mathbf I} \def\bfJ{\mathbf J}
 \def\bfK{\mathbf K} \def\bfL{\mathbf L} \def\bfM{\mathbf M} \def\bfN{\mathbf N} \def\bfO{\mathbf O} \def\bfP{\mathbf P}
 \def\bfQ{\mathbf Q} \def\bfR{\mathbf R} \def\bfS{\mathbf S} \def\bfT{\mathbf T} \def\bfU{\mathbf U} \def\bfV{\mathbf V}
 \def\bfW{\mathbf W} \def\bfX{\mathbf X} \def\bfY{\mathbf Y} \def\bfZ{\mathbf Z} \def\bfw{\mathbf w}
 \def\R {{\mathbb R }} \def\C {{\mathbb C }} \def\Z{{\mathbb Z}} \def\H{{\mathbb H}} \def\K{{\mathbb K}}
 \def\N{{\mathbb N}} \def\Q{{\mathbb Q}} \def\A{{\mathbb A}} \def\T{\mathbb T} \def\P{\mathbb P} \def\G{\mathbb G}
 \def\bbA{\mathbb A} \def\bbB{\mathbb B} \def\bbD{\mathbb D} \def\bbE{\mathbb E} \def\bbF{\mathbb F} \def\bbG{\mathbb G}
 \def\bbI{\mathbb I} \def\bbJ{\mathbb J} \def\bbK{\mathbb K} \def\bbL{\mathbb L} \def\bbM{\mathbb M} \def\bbN{\mathbb N} \def\bbO{\mathbb O}
 \def\bbP{\mathbb P} \def\bbQ{\mathbb Q} \def\bbS{\mathbb S} \def\bbT{\mathbb T} \def\bbU{\mathbb U} \def\bbV{\mathbb V}
 \def\bbW{\mathbb W} \def\bbX{\mathbb X} \def\bbY{\mathbb Y} \def\kappa{\varkappa} \def\epsilon{\varepsilon}
 \def\phi{\varphi} \def\le{\leqslant} \def\ge{\geqslant}

\def\UU{\bbU}
\def\Mat{\mathrm{Mat}}
\def\tto{\rightrightarrows}

\def\Gr{\mathrm{Gr}}

\def\graph{\mathrm{graph}}

\def\O{\mathrm{O}}

\def\la{\langle}
\def\ra{\rangle}

\def\B{\mathrm B}
\def\Int{\mathrm{Int}}
\def\LGr{\mathrm{LGr}}


\def\I{\mathbb I}
\def\M{\mathbb M}
\def\T{\mathbb T}
\def\S{\mathrm S}

\def\Lat{\mathrm{Lat}}
\def\LLat{\mathrm{LLat}} 
\def\Mod{\mathrm{Mod}}
\def\LMod{\mathrm{LMod}}
\def\Naz{\mathrm{Naz}}
\def\naz{\mathrm{naz}}
\def\bNaz{\mathbf{Naz}}
\def\AMod{\mathrm{AMod}}
\def\ALat{\mathrm{ALat}}
\def\MAT{\mathrm{MAT}}
\def\Mar{\mathrm{Mar}}

\def\Ver{\mathrm{Vert}}
\def\Bd{\mathrm{Bd}}
\def\We{\mathrm{We}}
\def\Heis{\mathrm{Heis}}
\def\Pol{\mathrm{Pol}}
\def\Ams{\mathrm{Ams}}
\def\Herm{\mathrm{Herm}}

\def\bbot{{\bot\!\!\!\bot}}

\def\shr{{\,{\ov{\ov {r\vphantom{h}}}}\,}}
\def\shg{{\,{\ov{\ov {g\vphantom{h}}}}\,}}
\def\shh{{\,{\ov{\ov {h\vphantom{h}}}}\,}}
\def\shp{{\,{\ov{\ov {p\vphantom{h}}}}\,}}

\def\ls{(\!(}
\def\rs{)\!)}

 \newcommand{\gr}{\mathop {\mathrm {gr}}\nolimits}

\begin{center}
 \Large \bf
 
 Algebras  of conjugacy classes
 in symmetric groups and checker triangulated surfaces
 
 \large \sc
 
 \bigskip
 
 Yu.A.Neretin%
 \footnote{Supported by the grants FWF, P28421,  P31591.}
 
\end{center}

{\small
	In 1999 V.~Ivanov and S.~Kerov observed that structure constants of algebras of conjugacy classes
	of symmetric groups $S_n$ admit a stabilization (in a non-obvious sense) as $n\to \infty$.
	We extend their construction to a  class of pairs  of groups $G\supset K$ and  
	algebras of conjugacy classes of $G$ with respect to $K$.
	In our basic example,
	  $G=S_n \times S_n$, $K$ is the diagonal
	 subgroup $S_n$. In this case we
	  get a geometric description of this algebra.}



\section{Introduction}

\COUNTERS

{\bf \punct Notation.}
1) Denote by $\# A$ the number of elements of a set $A$.
By $\coprod A_j$ we denote the disjoint union of sets $A_j$

\sm

2) Denote by $J_n$ the set 
$$J_n=\{1,2,\dots, n\} \subset\N.$$
Let $Y$ be a finite or countable set. Denote by $S(Y)$ the group of all finitely supported
permutations of $Y$. 
By $S_n$ we denote the group of permutations of $J_n$, by $S_\infty$ 
the group $S(\N)$.  We regard groups  $S_n$ as subgroups in $S_\infty$.

\sm

3) Let $G$ be a group, $K\subset G$  a subgroup.
Denote by $K\backslash G/K$ the set of double cosets  of $G$ with respect to $K$,
i.e., the quotient of $G$ with respect to the equivalence relation
$$
g\sim h_1 g h_2, \qquad\text{where $g\in G$, $h_1$, $h_2\in K$.}
$$
Denote by $G/\!\!/K$ the set of conjugacy classes of $G$ with respect to $K$,
i.e., the quotient of $G$ with respect to the equivalence relation
$$ 
g\sim hgh^{-1}, \qquad\text{where $g\in G$, $h\in H$.}
$$ 

Consider the group $G\times K$ and the subgroup $\wt K$ consisting
of elements $(h,h)\in G\times K$, where $h$ ranges in $K$. 
We have a canonical identification
$$
\wt K\backslash (G\times K)/\wt K\simeq G/\!\!/K.
$$
Indeed, let $(g,h)\in G\times K$. The corresponding double coset 
contains 
$$(g,h)(h^{-1},h^{-1})=(gh^{-1}, 1)$$
 and also all elements of the form
$$(r,r) (gh^{-1}, 1) (r^{-1},r^{-1})=\bigl( r(gh^{-1}) r^{-1}, 1\bigr).$$

For a group $G$ consider its multiples $G^{(m)}=G\times \dots \times G$
($m$ times), by $\diag (G)=\diag_m(G)\subset G^{(m)}$ we denote  the diagonal, i.e., the set
of all tuples $(g,\dots, g)\in G^{(m)}$. We have  obvious 
identifications
\begin{align*}
\diag_2(G)\backslash (G\times G)/ \diag_2(G)
&\simeq G/\!\!/G;
\\
\diag_3(G)\backslash (G\times G\times G)/ \diag_3(G)
&\simeq (G\times G)/\!\!/\diag_2(G),
\end{align*}
etc. 

\sm

4) For a finite group $G$ denote by $\C(G)$ the group algebra of
$G$, we denote the convolution by
 $*$. For $f\in \C(G)$ denote by $f^\star$ the function
$$
f^\star(g):=\ov{f(g^{-1})}.
$$
Clearly $f\to f^\star$ is an anti-involution of the group algebra
$$
(f_1*f_2)^\star= f_2^\star * f_1^\star, \qquad (f^\star)^\star=f.
$$

Denote by $\C(K\backslash G/K)$ (resp. $\C(G/\!\!/K)$)
the subalgebra of the group algebra consisting of 
functions that are constant on double cosets (resp. conjugacy classes),
these subalgebras are closed with respect to the anti-involution.


\sm

Let $\wt K\subset G\times K$ be as above. We have an obvious isomorphism
$$
\C (G/\!\!/K)\simeq \C\bigl(\wt K\backslash (G\times K)/ \wt K\bigr),
$$
i.e., convolution algebras of conjugacy classes are special cases of algebras
of double cosets.

\sm

{\bf \punct Bibliographical remarks on 
	 algebras of double cosets and their infinite-dimensional degenerations.}
 Formally, this subsection is not necessary, however here we explain origins and purposes
 of this work.
 
Let $\rho$ be a unitary representation of a finite group $G$ in a space $V$, denote by 
the same symbol $\rho$ the corresponding representation of the group algebra.
Denote by $V^K$ the subspace  of all $K$-fixed vectors, by $(V^K)^\bot$
its orthocomplement. The convolution algebra $\C(K\backslash G/K)$
acts in $V=V^K\oplus (V^K)^\bot$ by operators of block form
$$
\rho(f)=\begin{pmatrix}
\rho'(f)&0\\
0& 0
\end{pmatrix}
.$$
Thus for any representation $\rho$ of $G$ we have a representation
$\rho'$ of $\C(K\backslash G/K)$ in $V^K$. It can be easily shown that if $\rho$ is irreducible and $V^K\ne 0$,
then $\rho'$ 'remembers' $\rho$. For this reason convolution algebras
became tools of investigation of representations. We recall some well-known 
examples.

\sm

1) Hecke algebras (Iwahori \cite{Iwa}).
Let $G$ be the group $\GL(n,\mathbb{F}_q)$ of all 
invertible matrices over a finite field $\mathbb{F}_q$ and
$K$ be the group of upper-triangular matrices. These 
algebras  admit   explicit descriptions
and an interpolation in $q$ (for $q=1$ we get $\C (S_n)$). They can live
their  own lives independently of the group $\GL(n,\mathbb{F}_q)$, see, e.g., 
\cite{Bump}. For the subgroup $T$ of 
strictly upper triangular matrices the convolution algebra 
$\C\bigl(T\backslash \GL(n,\mathbb{F}_q)/T\bigr)$ also admits a transparent description,
see Yokonuma \cite{Yok}.

\sm

 There are some widely explored examples of convolution algebras related
to locally compact groups (in this case $K$ must be compact):

\sm

2) Let $G$ be a reductive Lie group and $K$  its maximal compact subgroup. 
Algebras $\C(K\backslash G/K)$ were widely  explored in classical representation theory of Lie groups
at least after Gelfand \cite{Gel}.
If $G$ is a rank one classical group, i.e., pseudounitary group
$\SU(1,n;\bbK)$, where $\bbK$ is $\R$, $\C$, or quaternions
(in usual notation,
$G=\SO(1,n)$, $\SU(1,n)$, $\Sp(1,n)$), then multiplication in this algebra
is determined by an explicit hypergeometric kernel, these algebras  have an explicit
two-parametric interpolation with respect to $n$ and a dimension $d=1$, $2$, $4$
 of a field
$\bbK$,  this interpolation
also includes one real form of the exceptional group $F_4$, see, 
Flensted-Jensen, Koornwinder \cite{FK}, see also \cite{Koo}.

\sm

3) Affine Hecke algebras (Iwahori, Matsumoto  \cite{IwaMats}). Let $\Q_p$
be a $p$-adic field, $\mathbb{O}_p$ be the ring of $p$-adic integers.
 Let $G$ be the group $\GL(n, \Q_p)$,
 and $K$ be the Iwahori subgroup. Recall that the Iwahori
 subgroup is a subgroup on $\GL(n,\mathbb{O}_p)$ consisting
 of matrices whose elements under the diagonal are contained
 in $p \mathbb{O}_p$. Such algebras $\C(K\backslash G/K)$ also admit an
 explicit description and an interpolation in
 $p$ and can live their own lives.

\sm

These examples have  further  extensions, however
in all these cases subgroups $K$ are quite large in $G$,
for smaller subgroups algebras of double cosets usually seem to be too
complicated.

\sm

It appears that for infinite-dimensional groups double coset
spaces $K\backslash G/K$  often admit a natural structure of a semigroup,
and for each unitary representation of $G$ this semigroup acts in the subspace
of $K$-fixed vectors. First example of such semigroup was discovered
by Ismagilov in \cite{Ism}. Many cases were examined by Olshansky \cite{Olsh-howe},
\cite{Olsh-symm},
he showed that  semigroups $K\backslash G/K$  admit explicit descriptions 
in some cases when their finite-dimensional counterparts seem 
non-handable. In \cite{Ner-book} it was observed that existence of such semigroup
structures 
 is a relatively usual phenomenon.  In \cite{Ner-char}--\cite{Ner-faa}, \cite{GN}  descriptions of such semigroups
were proposed on a quite general setting. 

\sm

A basic example
in \cite{Ner-umn} was semigroups $K\backslash G/K$, where
$$
G=S_\infty \times S_\infty \times S_\infty
$$
and $K$ is a subgroup in the diagonal $S_\infty$ fixing points 
$1$, $2$, \dots, $\alpha\in \N$.

The present paper is an attempt to look from infinity to finite objects.
We present a kind of  description (or quasi-description) of a family of  convolution semigroups
 related to symmetric groups, our basic example is
$$
\C\Bigl(\diag_3(S_n)\backslash (S_n\times S_n\times S_n)/\diag_3(S_n)\Bigr)
\,\,\simeq\,\, \C\Bigl((S_n\times S_n)/\!\!/\diag_2(S_n)\Bigr).
$$
We use arguments of Ivanov and Kerov \cite{IK} who observed that
the algebras $\C(S_n/\!\!/S_n)$ admit a  stabilization
(in a non-obvious sense) as $n\to\infty$. 

\sm 

{\bf \punct Ivanov-Kerov algebra of conjugacy classes of groups $S_j\times S_j$ with respect to the diagonal subgroup.}
Denote
$$
G_j:=S_j\times S_j, \qquad K_j:=\diag_2(S_j).
$$
For
$$
\text{$g$, $h$, $r$, \dots $\in \coprod_{j=0}^\infty G_j$}
$$
denote by
  $$\text{$\shg$, $\shh$, $\shr$, \dots $\in \coprod_{j=0}^\infty G_j/\!\!/K_j$}
$$
the corresponding conjugacy classes.

Let $j\le N$ and $g\in G_j$. Denote by $\wt g$ the corresponding
element of $G_N\supset G_j$.  We define an element $A_N[\shg]$
of $\C(G_N/\!\!/K_N)$ by the formula 
$$
A_N[\shg]=A_N^{(j)}[\shg]:=\frac 1{(N-j)!} \sum_{\tau \in K_N} 
\tau^{-1} \wt g \tau.
$$
Sometimes we write superscript $^{(j)}$ to emphasize that $\shg$ is considered as an element of $G_j/\!\!/ K_j$.
If $j>N$ and $g\in G_j$, we set
$$
A_N[\shg]:=0.
$$

Thus for each $N$ we get a family of elements of $\C (G_N/\!\!/K_N)$.

\sm

{\sc Remark.} 
Elements $A_N[\shg]=A_N^{(N)}[\shg]$, where $\shg$ ranges in $G_N/\!\!/K_N$, form a basis
in $\C (G_N/\!\!/G_N)$. However, if  $g\in G_N$ actually
is an element of a subgroup $G_k$, our family contains also elements
\begin{multline*}
A_N^{(N-1)}[\shg]=\frac 1{1!}A_N[\shg],\quad A_N^{(N-2)}[\shg]=\frac 1{2!}A_N[\shg],\quad\dots,
\\
 A_N^{(N-k)}[\shg]=\frac 1{(N-k)!}A_N[\shg].
\end{multline*}


\begin{theorem}
	\label{th:1}
	 Let $\shg$, $\shh$, $\shr$ range in the disjoint union $\coprod_{j=0}^\infty G_j/\!\!/K_j$.
	 Then there are 
	non-negative  integers $a_{\shg,\shh}^\shr$, which do not depend on $N$, satisfying the following properties:
	
	\sm
	
$\bullet$   For each $N$
	\begin{equation}
	A_N[\shg]*A_N[\shh]=\sum_\shr a_{\shg,\shh}^\shr A_N[\shr].
	\label{eq:0}
	\end{equation}
	
	\sm
	
	$\bullet$ Consider a linear space $\cB$ with a basis consisting of symbols $A[\shg]$,
	where $\shg$ ranges in $\coprod_{j=0}^\infty G_j/\!\!/K_j$. Then the
	formula 
	\begin{equation}
	A[\shg]*A[\shh]=\sum_\shr a_{\shg,\shh}^{\shr} A[\shr],
	\label{eq:1}
	\end{equation}
	determines a structure of an associative algebra on $\cB$.
\end{theorem}



\begin{figure}
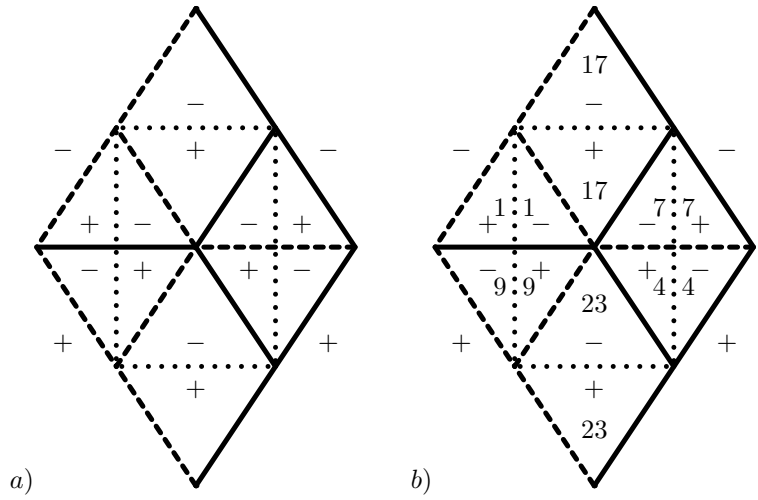

	$$a)\epsfbox{kerov-ris1.3}\qquad b) \epsfbox{kerov-ris1.4}$$
	\caption{
		a) A piece of a checker triangulated surface.
		\newline
	b)	A piece of a labeled checker triangulated surface.
		\label{ris:3}}	
\end{figure}

\begin{figure}
	$$\epsfbox{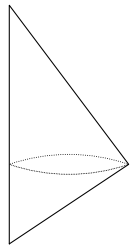}$$
	\caption{A surface consisting of two triangles.%
	\label{fig:triangle}}
\end{figure}	

\begin{figure}
	$$\epsfbox{kerov-ris1.5}$$
	\caption{Surfaces $\cR$, $\cQ$ and a partial bijection $\cR_-$ to $\cQ_+$.\label{ris:4}}
\end{figure}
%

In Theorem \ref{th:1-1} we present a geometric description (or quasi-description) of this algebra.

\sm

{\bf\punct Checker triangulated surfaces and symmetric groups.}
We say that a {\it checker triangulated surface} is a 
 two-dimensional oriented closed (generally, disconnected) surface 
 (see Fig. \ref{ris:3}) equipped with
 the following data:
 
 \sm
 
$\bullet$ a graph separating
 surface into triangles ({\it faces});

\sm

$\bullet$ for each triangle we assign a sign $(+)$ or $(-)$, pluses and minuses
are arranged  in a checker order;

\sm

$\bullet$ edges are colored by red, yellow, blue; colors of edges of each face are different;
in  plus-triangles these colors are located clockwise, on minus-triangles anti-clockwise.

\sm

We admit non-connected surfaces.  Also, we admit  pairs of triangles, which
have two or three common edges as on Fig. \ref{fig:triangle}. 

Surfaces are defined up to the combinatorial equivalence.
Denote by $\Xi_N$ the set of all checker triangulated surfaces with $2N$ faces.

We say that a {\it labeling} of a surface $\cR\in \Xi_N$ is a bijective map
from the set $J_N=\{1,2, \dots, N\}$ to the set of all plus-triangles. We automatically
assign labels to minus-triangles assuming that triangles separated by {\it blue}
edges have same labels. Denote by $\wt \Xi_n$ the set of all labeled surfaces with $2N$
triangles defined up to a combinatorial equivalence%
\footnote{If a surface $\cR$ admits combinatorial automorphisms, then different labelings of $\cR$ can give
equivalent labeled surfaces.}.

\sm

{\it There is a natural one-to-one correspondence between the set 
	 $\wt\Xi_N$ and the group $G_N=S_N\times S_N$.}

\sm

Indeed, fix a  labeled   checker surface $\wt \cR$. For each red edge $v$ we
consider labels $i_+(v)$ and $i_-(v)$ on its plus and minus sides.
Assuming $\sigma_{red}:i_+(v)\mapsto i_-(v)$ we get an element $\sigma_{red}$ of $S_N$.
Considering yellow edges we obtain another permutation $\sigma_{yellow}\in S_N$.

A permutation $\tau\in S_N$ of labels corresponds to a simultaneous conjugation
$$
(\sigma_{red}, \sigma_{yellow})\mapsto 
(\tau\sigma_{red} \tau^{-1}, \tau \sigma_{yellow}\tau^{-1}).
$$
Therefore
{\it we get a canonical bijection between sets $\Xi_N$ and $(S_N\times S_N)/\!\!/S_N$.}

\sm

\begin{figure}
	\epsfbox{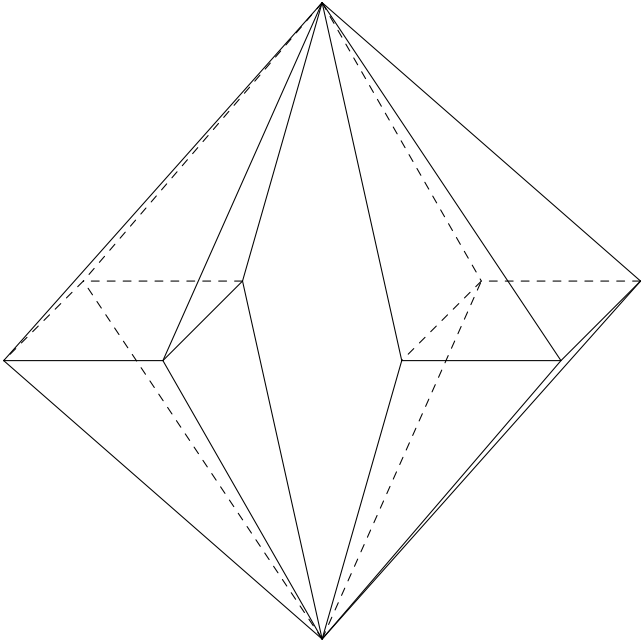}
	\caption{A surface whose vertices are glued.\label{ris:5}}
\end{figure}

Next, let us describe a product in $S_N\times S_N$ in terms of labeled surfaces.
Let $\wt\cP$ and $\wt \cR$ be labeled surfaces. For each $j=1$, \dots, $N$,
we identify the $j$-th minus-triangle of $\wt R$ with the $j$-th plus-triangle
of $\cP$ according colors of their sides. Thus we get a two-dimensional
simplicial cell complex consisting of $N$ labeled plus-triangles inherited from
$\cR$, labeled $N$ minus-triangles inherited from $\cP$ and $N$ {\it plus-minus-triangles}
obtained as result of gluing. Each edge is contained in 3 triangles.
Removing interiors of 
all plus-minus-triangles we come to a simplicial cell complex such that
each edge is contained in two triangles. 
 In fact, this is a surface, but some vertices of the surface can be glued
one with another (see Fig. \ref{ris:5}). Cutting all such gluings we get a new  surface whose triangles are equipped with labels. Some of blue edges were contained
in removed triangles. Labels on both sides of such edge coincide with a label on the
removed triangle and therefore coincide. So  our new surface is correctly labeled. 

\sm

{\bf \punct Remark. Belyi data.} Checker triangulated surfaces arise in a natural
way in algebraic geometry under names {\it Belyi data} or
{\it dessigns d'enfant}.

Denote by $\ov \Q$ the field of algebraic numbers.
Consider an algebraic curve $C$ and a meromorphic function ({\it Belyi function})
$f$ on $C$ whose critical values are contained in the set $0$, $1$, $\infty$. 
According the famous Belyi theorem \cite{Bel}  such a function on a given curve $C$
exists if and only if $C$ can be determined by
a system of algebraic equations with coefficients in $\ov\Q$.

Consider the Riemann sphere $\ov\C=\C \cup \infty$
and the real projective line $\ov\R=\R\cup\infty$ in $\ov \C$.
Let us say that the upper half-plane is a plus-triangle, lower half-plane is a minus-triangle, the segment
$[1,\infty]$ is red, the segment $[0,1]$ is yellow,
and the segment $[-\infty,0]$ is blue. Thus the Riemann
sphere $\ov\C$ becomes a checker triangulated surface.
The preimage of $\R$ on $C$ is colored graph splitting $C$
into triangles. Clearly, we come to a checker triangulated surface. 

The Galois group of $\ov \Q$ over $\Q$ acts on the set of all Belyi functions.
A.~Gro\-then\-dieck proposed a program of investigation of the Galois groups using 
Belyi functions and graphs on surfaces, see,
 e.g.,  \cite{ShV}, \cite{Sch}, \cite{S-Loshak}, \cite{G-G}, \cite{Shabat}.

 Relations of this topic and infinite symmetric
 group remain to be non-clear.
 
 \sm
 
 {\bf\punct Partial bijections.}
 Recall that a {\it partial bijection} $\lambda$ of a set $Y$ to a set $Z$ is a bijection of a subset
 $A\subset Y$ to a subset $B\subset Y$ (we admit the case $A=B=\varnothing$). We define  rank, domain, and image of a partial bijection by
 $$
 \rk \lambda:= \# A=\#B, \qquad \dom \lambda:=A, \qquad \im \lambda:=B.
 $$
 Denote by $\PB(Y,Z)$ the set of all partial bijections $Y\to Z$.

 For partial bijections $\mu:W\to Y$, $\lambda: Y\to Z$ we define their product $\lambda\mu:W\to Z$.
 We say that $w\in \dom \lambda\mu$ if $w\in \dom \mu$ and $\mu (w)\in \dom \lambda$.
 In this case we set $\lambda\mu (w)=\lambda(\mu (w))$
 
 \sm

 {\bf \punct Construction  of the algebra $\cB$.}
 The basis of the algebra is numerated by the set
 $$
 \coprod_{n=0}^\infty \Xi_n.
 $$
Let $\cR\in \Xi_n$, $\cQ\in \Xi_k$.  Let $\lambda$ be a partial bijection from the set
of minus-triangles of
$\cR$ to the set of plus-triangles of $\cQ$  (see Fig. \ref{ris:4}).
Consider the disjoint union $\cR\coprod \cQ$ and let us perform the following transformations.
For each face $A\in \dom \lambda$, we take the face $\lambda(A)\in \im\lambda$,
remove both faces and identify their boundaries according colors of edges. 
After this, we get a compact two-dimensional simplicial complex, and each edge of the complex is contained
in precisely two faces.
As above
  some vertices of the surface can be glued
one with another (see Fig. \ref{ris:5}). Cutting all such gluings we get a new surface 
$$
\cR \circledast_\lambda\cQ.
$$

In this notation, the product is defined by
$$
\cR *\cQ=\sum_{\lambda\in \PB(\cR_-,\cQ_+)} \cR \circledast_\lambda\cQ.
$$ 

\begin{theorem}
\label{th:1-1}
	This algebra coincides with the algebra $\cB$ 
	defined in Theorem {\rm \ref{th:1}}.
\end{theorem}

\section{A more general construction}

\COUNTERS

\begin{figure}
	$$\epsfbox{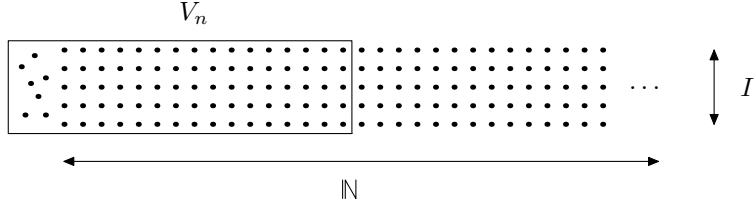}$$
	\caption{Sets $V\supset V_n$.\label{ris:1}}
\end{figure}

{\bf\punct Groups $G_n$.} Fix a finite set $I$.
Consider a finite or countable set%
\footnote{The case $X=\varnothing$ is sufficiently interesting.} $X$, the product $I\times \N$
and the unions
$$
V=X\cup (I\times \N),\qquad V_n=X \cup (I\times J_n)\subset V,
$$
see Fig. \ref{ris:1}.
The group $S_\infty$ acts on $\N$, therefore it acts on the product $I\times \N$.
We can imagine $I\times \N$ as a table with $\#I$ infinite rows.
The group 
$S_\infty$ acts by permutations of columns. 

Next, we consider the trivial action of $S_\infty$ on $X$,
this determines an action  of $S_\infty$ on the whole  $V$.

Consider a subgroup $G_\infty\subset S(V)$ containing $S_\infty$.
Set
$$
G_n:= G_\infty\cap S(V_n).
$$
In other words, we consider elements of $G$ supported by $V_n$.
As above, we denote
$$K_n:=S_n, \qquad K_\infty=S_\infty.$$

\sm

{\bf \punct Some natural examples.}
a) The case discussed above corresponds to $X=\varnothing$ and a 2-element set $I$.
The group $G_\infty=S_\infty\times S_\infty$ acts by permutations preserving rows.

\sm

b) $X=\varnothing$, $G=S(I\times \N)$.

\sm 

c) Let $X=\varnothing$, $G=S_\infty\times \dots\times S_\infty$
consists of permutations of $I\times V$ preserving each row.

\sm

d) $X$ is a finite set, $\#I=1$, $G=S(X\cup \N)$.

\sm

e) $X=\varnothing$, $\# I=k$, 
 $G$ is a semidirect product of $S_\infty$ acting by permutations
of columns and the group $(S_k)^\infty$ acting by permutations of elements
of each column.  

\sm

In all these cases the sets of conjugacy classes admit geometric descriptions
by tricks of \cite{Ner-umn}.

\sm

{\bf \punct Algebras $\C(G_n/\!\!/K_n)$.}  Let $\shg$, $\shh$, $\shr$, \dots range 
in $\coprod_{j=0}^\infty G_j/\!\!/K_j$. 
Denote by $g$, $h$, $r$ their representatives.

Fix $n=0$, $1$, $2$, \dots, let $\shg\in G_n/\!\!/K_n$, $g\in\shg$, and   $N\ge n$.
For any $n$-element set $\Omega\subset \N$ we define the sum
$$
R(g,\Omega):=\sum_{\sigma:J_n\to\Omega} \sigma g \sigma^{-1}\quad\in \C(G_\infty)
,$$
where $\sigma$ ranges in the set of bijective maps $J_n\to \Omega$.
We explain this more carefully.
Any $\sigma$ determines a bijection $I\times J_n\to I\times \Omega$,
we denote it by the same symbol $\sigma$. Thus, 
$ \sigma g \sigma^{-1}$ is a map
$$
\sigma g \sigma^{-1}:\, X\cup (I\times \Omega) \to X\cup (I\times \Omega).
$$
We extend it to a map $V_N\to V_N$ in a trivial way,
$$
\sigma g \sigma^{-1} (w):=w,\qquad \text{if $w\in I\times (J_N\setminus\Omega$).}
$$

Next, we consider the sum
$$
A_N[\shg]=\sum_{\Omega\subset J_N:\,\, \#\Omega=n} R(g,\Omega).
$$

For $N< n$, we set $A_N[\shg]=0$.

\sm

{\sc Remark.} Clearly, for $N\ge n$, we have
\begin{equation}
A_N[\shg]=\frac 1{(N-n)!} \sum_{\tau\in K_N} \tau g\tau^{-1}.
\label{eq:A-con}
\end{equation}
However, our long  definition will be used below.
\hfill $\boxtimes$

\begin{theorem}
	\label{th:2}
	For the groups $G_n$ defined in this subsection Theorem {\rm\ref{th:1}}
	remains true.
\end{theorem}

Thus, for any group $G_\infty$, we get an associative algebra $\cB=\cB[G_\infty/\!\!/K_\infty]$
and canonical epimorphisms
$\cB[G_\infty/\!\!/K_\infty]\to \C(G_N/\!\!/K_N)$.

\sm

{\bf\punct Local bijections.} Following Ivanov and Kerov \cite{IK} we define a semigroup
of local bijections%
\footnote{See \cite{AN}, \cite{MMN}, \cite{Mel}, \cite{KR} on continuations of the work \cite{IK}.}.
Let $Y$ be a finite or countable set.
A {\it local  bijection} $\omega$ is a bijection from a
finite subset $\Omega\subset Y$ to $\Omega$. We denote such local bijection by
$\ls\omega, \Omega\rs$. Any local bijection admits a canonical extension to an
element $\wt\omega\in S(Y)$, we set $\wt\omega y=y$ if $y\notin \Omega$.
We define the product of local bijections by
$$
\ls\omega,\Omega\rs\circ \ls\mu,M\rs=\Bigl(\!\!\Bigl(\wt \omega \cdot \wt\mu \Bigr|_{\Omega\cup M},
\Omega\cup M\Bigr)\!\!\Bigr).
$$
Denote by $\cL(Y)$ the semigroup of all local bijections of $Y$.

\sm

{\sc Remark.} A local bijection determines a partial bijection. But the $\circ$-product 
differs from 
the product of partial bijections.
\hfill $\boxtimes$

\sm 

The group $S(Y)$ acts on $\cL(Y)$ by conjugations in the obvious way,
$$
\sigma\ls w,\Omega \rs \sigma^{-1}=\ls \sigma \wt w \sigma^{-1}\Bigr|_{\sigma\Omega}, \sigma\Omega  \rs.
$$

On the other hand, we have a natural forgetting homomorphism 
$$\iota:\cL(Y)\to S(Y)$$
defined by 
$$
\iota\ls\omega,\Omega\rs=\wt\omega.
$$

{\bf\punct The semigroup algebra for local bijections.}
Denote by $\C[\cL(Y)]$ the space of all formal series of the form
$$
\sum_{\ls\omega,\Omega\rs\in \cL(Y)} a_{\ls\omega,\Omega\rs} \ls\omega,\Omega\rs,\qquad a_{\ls\omega,\Omega\rs}\in\C.
$$
This space is equipped with a convolution
$$
\sum_{\ls\omega,\Omega\rs} a_{\ls\omega,\Omega\rs}\, \ls\omega,\Omega\rs \, * 
\sum_{\ls\mu,M\rs} b_{\ls\mu,M\rs}\, \ls\mu,M\rs
= \sum_{\ls\nu,N\rs} c_{\ls\nu,N\rs}\, \ls\nu,N\rs,
$$
where
$$
c_{\ls\nu,N\rs}=\sum_{\ls\omega\Omega\rs, \ls\mu M\rs: \,\, \Omega\cup M=N,\,\,
	\wt\omega \wt \mu=\wt \nu} a_{\ls\omega,\Omega\rs} b_{\ls\mu,M\rs}
$$
(this sum is finite).

For any subset $Z$ in $Y$ we have a homomorphism 
$$\pi^Y_Z:\C[\cL(Y)]\to \C[\cL(Z)]$$
defined
on generators by
$$
\pi\ls\omega,\Omega\rs=\begin{cases}
\ls\omega,\Omega\rs,\qquad \text{if $\Omega\subset Z$};
\\
0,\qquad \text{otherwise}.
\end{cases}
$$

For a {\it finite} set $Y$ 
we have a homomorphism of algebras 
$$
\iota: \C[\cL(Y)]\to \C[S(Y)]
$$
defined by
$$
\iota\Bigl( \sum_{\ls\omega,\Omega\rs} a_{\ls\omega,\Omega\rs}\, \ls\omega,\Omega\rs\Bigr)=\!\!
\sum_{\ls\omega,\Omega\rs} a_{\ls\omega,\Omega\rs} \wt\omega
=\!\!\sum_{g\in S(Y)} \Bigl[\sum_{\Omega:\, \Omega\supset \{\text{support of $g$} \}}
a_{\ls g\Bigr|_\Omega,\Omega \rs} \Bigr]\,g
$$
(for an infinite set $Y$ the sum in square brackets  is  infinite).

\sm





{\bf \punct Elements $B[\shg]$.} Let $\shg\in G_n/\!\!/K_n$, $g\in \shg$. For any   $n$-element subset $\Omega\subset \N$
we define an element of $\C[\cL(V)]$ by
$$
R(\shg,\Omega) =
\sum_{\sigma:J_n\to \Omega} \bigl(\!\bigr(\sigma g\sigma^{-1}, X\cup (I\times \Omega)\,\bigr)\!\bigr)
,
$$
where the summation is taken over all bijections $\sigma:J_n\to \Omega$.
Equivalently, we can chose one bijection $\sigma_0:J_n\to \Omega$
and write the formula in the form
$$
R(\shg,\Omega) =
\sum_{u\in S_n} \bigl(\!\bigr(\sigma_0 u g u^{-1}\sigma_0^{-1},
 X\cup (I\times \Omega)\,\bigr)\!\bigr).
$$

Next, we define
$$
B[\shg]:=\sum_{\Omega:\,\#\Omega=n}  R(\shg,\Omega).
$$
By construction these elements are invariant with respect to conjugations by elements of $K_\infty$,

\sm

{\sc Example.} 
Let $G_\infty$ be $S_\infty \times S_\infty$ acting on $\{1,2\}\times \N$ as above. For  $\shg\in (S_N \times S_N)/\!\!/S_N$
we take the corresponding
  $\cR\in \Xi_N$. Then
 $$
 B[\shg]=\sum_{\Omega\subset\N:\, \#\Omega=N} \sum_{\kappa} \wt\cR_\kappa,
 $$
  where the summation in the interior sum is taken over all bijective maps $\kappa$ from $\Omega$
  to the set of plus-triangles of $\cR$. A checker triangulated surface equipped with such a map
  determines an element of the group $S(\Omega)\times S(\Omega)$. We consider this element as a local
  bijection of $\{1,2\}\times \N$ with domain $\{1,2\}\times \Omega$.
  \hfill $\boxtimes$

\sm

For any $\shg\in G_n/\!\!/K_n$, $\shh\in G_m/\!\!/K_m$ the convolution $B[\shg] * B[\shh]$
is $K_\infty$-invariant, therefore it has the form
$$
B[\shg]*B[\shh]:=\sum_{\shr\in \coprod G_j/\!\!/K_j} a^{\shr}_{\shg,\shh} B[\shr].
$$

{\bf\punct Proof of Theorem \ref{th:2}.}
First, we apply the map $\pi^V_{V_N}$.
For any $\shg$ and  $n\in \N$ we get an element
$$
B_N[\shg]:=\pi^V_{V_N} B[\shg].
$$
Since $\pi^V_{V_N}$ is a homomorphism of algebras, we get
$$
B_N[\shg]*B_N[\shh]:=\sum_{r\in \coprod_{k=0}^N G_j/\!\!/K_j} a^{\shr}_{\shg,\shh} B_N[\shr], 
$$
where $g$, $h$ range in $\coprod_{k=0}^N G_k/\!\!/S_k$. 

Next, we apply the forgetting map $\iota:\C[\cL(V_n)]\to\C[S(V_n)] $ to elements 
$B_N[\shg]$. Evidently, we get
$$
\iota (B_N[\shg])=A_N[\shg].
$$
This implies our statement.

\sm

{\bf \punct An expression for  products.}
Let $g\in G_n$, $h\in G_k$. 
Let $\lambda$ ranges in the set of  partial bijections $J_k\to J_n$. Fix $\lambda$, denote $d=\rk \lambda$.
Fix any pair of  injective maps $\sigma_0: J_n\to J_{n+k-d}$, $\tau_0:J_k\to  J_{n+k-d}$
such that the product $\sigma^{-1}_0\tau_0$ as a product of partial bijections 
is $\lambda$.
Define the conjugacy class 
$$
\ov{\ov{g\circledast_\lambda h}}:= 
\ov{\ov{\sigma_0 \ls g,V_n\rs \sigma^{-1}_0\circ \tau_0 \ls h,V_k\rs\tau^{-1}_0 }}\,\in G_{n+k-d}/\!\!/S_{n+k-d},
$$
it
does not depend on the choice of $\sigma_0$ and $\tau_0$.

\begin{theorem}
\label{th:3}
	\begin{equation}
	B[\shg]*B[\shh]=\sum_{\lambda\in \PB(J_k,J_n)}  B[\ov{\ov{g\circledast_\lambda h}}].
	\label{eq:BBB}
	\end{equation}	 
\end{theorem}

\begin{figure}
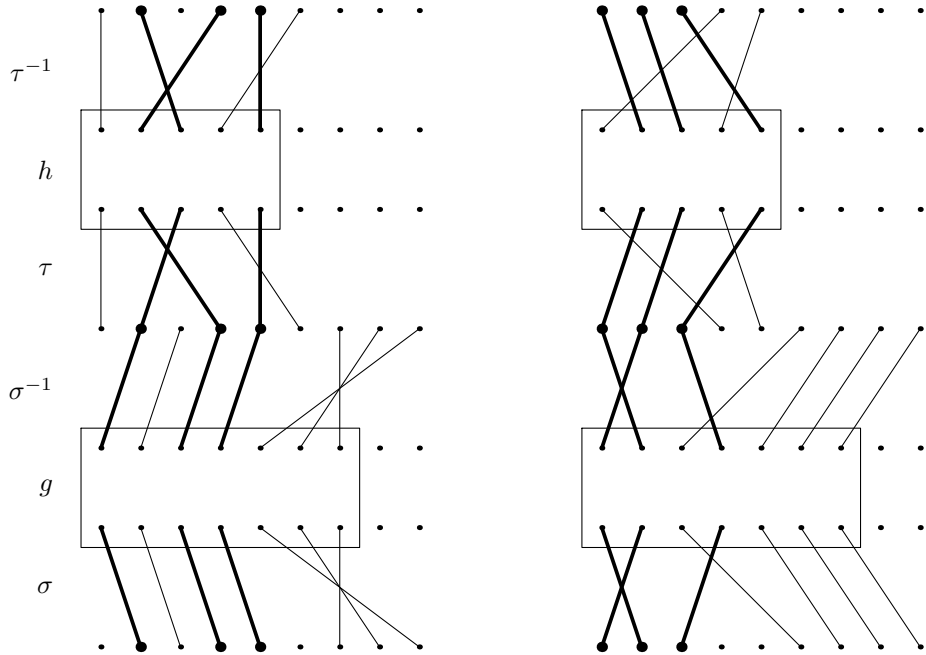

	\epsfbox{kerov-ris1.1}\qquad \qquad\qquad \epsfbox{kerov-ris1.2}
	
	\caption{The left part. A product $\sigma \ls g,V_n\rs \, \sigma^{-1}\circ \tau \ls h,V_k\rs\, \tau^{-1}$.
		Here $k=5$, $n=7$, $d=3$, $n+k-d=9$.
		We draw 7 copies of the set $J_{n-k+1}$. Maps $\tau^{-1}$, $\tau$, $\sigma^{-1}$, $\sigma$
		are shown by arcs. We mark  $\sigma \dom\lambda=\tau^{-1}\im\lambda$
		as fat points. Arcs containing fat points are thick. Boxes correspond to elements  
		$\ls h, V_k\rs\in G_k$ and $\ls g, V_n\rs\in G_n$. 
		\newline
		The right part. A canonical form of a pair $(\tau,\sigma)$.
		\label{ris:2}
	}
\end{figure}

{\sc Proof.} We expand the product $B[\shg]*B[\shh]$ by definition as a sum
\begin{equation}
\sum_{\sigma,\tau} \sigma \ls g,V_n\rs \sigma^{-1}\circ \tau \ls h,V_k\rs \tau^{-1}=
\sum_{\sigma, \tau}  \ls \sigma g\sigma^{-1}, \sigma V_n\rs
\circ
\ls \tau h \tau^{-1}, \tau V_k\rs,
\label{eq:ts}
\end{equation}
where $\sigma:J_n\to \N$, $\tau:J_k\to \N$ are injective maps. We wish to show that this expansion
coincides with (\ref{eq:BBB}). It is sufficient to identify the following sub-sum of (\ref{eq:ts})
\begin{equation}
\sum_{\sigma, \tau: \,\,\sigma J_n\cup \tau J_k= J_{n+k-d}}  \ls \sigma g\sigma^{-1}, \sigma V_n\rs 
\circ  \ls \tau h \tau^{-1}, \tau V_k \rs
\label{eq:subsum}
\end{equation}
and the sub-sum of (\ref{eq:BBB}) consisting of all summands of the form 
\begin{multline}
\sum_{\lambda\in \PB(J_k,J_n )} R (\ov{\ov{g\circledast_\lambda h}},J_{n+k-d})
=\\= \sum_{\lambda\in \PB(J_k,J_n )} \sum_{u\in S_{n+k-d}} 
\ls u^{-1}(g\circledast_\lambda h)u, V_{n+k-d}\rs.
\label{eq:subsum2}
\end{multline}

Notice that the group $S_{n+k-d}$ acts on the set of summands of (\ref{eq:subsum2})
permuting summands of the interior sum. The corresponding
action on (\ref{eq:subsum})
has the form
$$
u:\,\, (\sigma,\tau)\mapsto (u^{-1}\sigma ,u^{-1}\tau ).
$$
We refer to Fig. \ref{ris:2}. The second action corresponds to simultaneous application 
of a substitution $u$ to rows number 1, 4, 7. Rows 2, 3, 5, 6 remain to be fixed, and arcs are  moved by
corresponding permutations of their ends.
Applying an appropriate $u$ we can put fat points to positions 1, 2, \dots, $d$. Moreover, we can make $\tau^{-1}u$
monotone on $J_d$:
$$
i<j\le d\qquad \Rightarrow \qquad \tau^{-1}u (i)< \tau^{-1}u (j).
$$ 
(on Fig. \ref{ris:2} this means that the corresponding arcs have no intersections).

Next, we can put points of $\tau J_k\setminus \tau J_d$ to points of $J_k\setminus J_d$.
Moreover, we can make $\tau^{-1}u$ monotone on $J_k\setminus J_d$.
Finally, we can  put points of $\sigma J_n\setminus \sigma J_d$ to points of $J_{n+k-d}\setminus J_k$.
Moreover, we can make $\sigma^{-1} u$ monotone on $J_{n+k-d}\setminus J_k$.

This determines $u$ in a unique way. On the other hand, the partial bijection $\sigma^{-1} \tau$
does not changed under this transformation, and the new pair $(u^{-1}\sigma ,u^{-1}\tau)$ is uniquely
determined by $\dom \sigma^{-1} \tau$, $\im \sigma^{-1} \tau$ and the map
$\sigma^{-1} u: J_d \to \im \sigma^{-1} \tau$. 

Thus, we see that orbits of the group $S_{n+k-d}$ on {\it the set of summands} 
of (\ref{eq:subsum}) are enumerated by partial bijections $\lambda$ and stabilizers are trivial.
\hfill $\square$ 


\section{Final remarks}

\COUNTERS

{\bf \punct The involution.} The map $g\mapsto g^{-1}$ determines an anti-involution
in the algebra $\cL[V]$ and anti-involution in $\cB[G_\infty/\!\!/K_\infty]$,
$$
B[\shg]^\star:=B[\ov{\ov{g^{-1}}}].
$$
Evidently,
$$
(B[\shg]* B[\shh])^\star=B[\shh]^\star * B[{\shg}]^\star.
$$

{\bf \punct The filtration.}
Fix $G_\infty$ and set
$$
\cB:=\cB[G_\infty/\!\!/K_\infty].
$$
Let $\cB_n$
be the subspace in $\cB$ generated by all
$$B[\shg], \qquad\text{where $g$ ranges in $\coprod_{j=0}^n G_j/\!\!/S_j$.} 
$$
We get an increasing filtration,
$$\dots\subset \cB_n\subset \cB_{n+1}\subset\dots$$
Evidently,
$$
V\in \cB_k, \quad W\in \cB_l\quad \Rightarrow\quad V*W\in \cB_{k+l}.
$$

{\bf \punct The associated graded algebra.}
We construct the graded algebra $\gr\cB[G_\infty/\!\!/K_\infty]$ in the usual way.
Namely, the product
$$
\cB_k\times \cB_l \,\, \to \cB_{k+l}
$$
determines a map
$$
\cB_k/\cB_{k-1}\,\times \cB_l/\cB_{l-1} \,\, \to \,\,\cB_{k+l}/\cB_{k+l-1}
.
$$
In this way, we get a structure of an associative algebra on
$$
\gr\cB[G_\infty/\!\!/K_\infty]=\bigoplus_{k=0}^\infty \cB_k/\cB_{k-1}.
$$
Each subspace $ \cB_k/\cB_{k-1}$
has a natural  basis enumerated by elements of $ G_k/\!\!/S_k$. 

It is easy to describe the multiplication  in $\gr\cB[G_\infty/\!\!/K_\infty]$: in the sum in the right-hand side
of (\ref{eq:BBB}) we leave only the first summand (corresponding to the partial bijection of rank 0).

Formulate this more precisely. For $n$, $k\in \Z_+$ denote by $\theta_{n,k}$ the partial bijection
$J_n\to J_{n+k}\setminus J_k$ defined by $j\mapsto j+k$.

\begin{proposition}
	{\rm a)} For $\shg\in G_n/\!\!/S_n$, $\shh\in G_k/\!\!/S_k$, their product in
	graded algebra $\gr\cB[G_\infty/\!\!/K_\infty]$
	is
	$$
	B[\shg]\diamond B[\shh]= B[\ov{\ov{\theta_{n,k} \ls g,V_n\rs \theta_{n,k}^{-1} \circ \ls h,V_k\rs}}]
	.
	$$
	
	{\rm b)} If $X=\varnothing$, then the algebra $\gr\cB[G_\infty/\!\!/K_\infty]$ is commutative.
\end{proposition}

In fact, we get a semigroup
structure on $\coprod_{j=0}^\infty G_j/\!\!/S_j$
given by
$$
\shg\bullet \shh:= \ov{\ov{\theta_{n,k} \ls g,V_n\rs \theta_{n,k}^{-1}\circ  \ls h,V_k\rs}}.
$$
The algebra $\gr\cB[G_\infty/\!\!/K_\infty]$ is the semigroup algebra of this semigroup.

\sm

{\sc Example.} For $G_\infty=S_\infty\times S_\infty$ the $\bullet$-product
corresponds to the disjoint union of checker triangulated surfaces.
\hfill $\boxtimes$

\sm

{\sc Remark.}
These semigroups are similar to semigroups of double cosets, which were considered
in \cite{Ner-umn}, \cite{Ner-imrn}.
However, degree of generality in \cite{Ner-umn} is wider
(as we noticed in Introduction,
conjugacy classes are special cases of double cosets and not vice versa),  even 
in the same situations we get
 slightly different structures. For instance, for $G_\infty=S_\infty\times S_\infty$
adding to a checker triangulated surface $\frR\in \Xi_j$ a collection of $k$ double triangles
(drawn on Fig. \ref{fig:triangle}) we get different objects (this corresponds to embeddings
$S_j\times S_j\to S_{k+j}\times S_{k+j}$). However, such elements in  \cite{Ner-imrn}
are identified. Notice that we have a natural identification of sets
$$
\coprod_{j=0}^\infty G_j/\!\!/ K_j\,\simeq (G_\infty/\!\!/K_\infty)\times \Z_+,
$$
where $\Z_+$ denotes the set of non-negative integers. The left-hand side enumerates elements of the basis
in $\cB[G_\infty/\!\!/K_\infty]$. The semigroup
$$
G_\infty/\!\!/K_\infty\,\simeq\, K_\infty \backslash G_\infty /K_\infty
$$
is one of objects of 
\cite{Ner-imrn}.  
\hfill $\boxtimes$

\sm

{\bf \punct The Poisson bracket.} Let $X=\varnothing$,
in particular, in this case $\gr \cB[G_\infty/\!\!/K_\infty]$ is commutative. Consider the map
$$
\cB_k \times \cB_l\,\, \to\,\, \cB_{k+l-1}
$$
given by
$$
(V,W)\mapsto [V,W]=V*W-W*V
.$$
As usual, we get a map 
$$
\cB_k/\cB_{k-1} \times \cB_l/\cB_{l-1} \,\, \to\,\, \cB_{k+l-1}/\cB_{k+l-2} 
$$
and a structure of a Lie algebra on the space $\gr\cB[G_\infty/\!\!/K_\infty]$.
It is easy to write a formula for the bracket 
$$
\bigl[B[\shg], B[\shh]\bigr]_{gr}=
\sum_{\lambda\in \PB(J_k,J_n):\, \rk\lambda=1}\Bigl( B[\ov{\ov{g\circledast_\lambda h}}]-
B[\ov{\ov{h\circledast_{\lambda^{-1}} g}}]\Bigr).
$$
Of course, a partial bijection $J_k\to J_n$ of rank 1 is
determined by a pair $\alpha\in J_k$, $\beta\in J_n$.

\sm

{\sc Remark.}
Recall that there is a well-known Poisson structure on spaces
$
(K\times \dots \times K)/\!\!/\diag{K}
$, where $K$ is a compact Lie group, see \cite{Gold}, \cite{FR}. 
\hfill $\boxtimes$

\noindent
\tt Math.Dept., University of Vienna,
\\
Oskar-Morgenstern-Platz 1, 1090 Wien;
\\
Institute for Information Transmission Problems;
\\
Institute for Theoretical and Experimental Physics (until 11.2021);
\\
Mech.Math.Dept., Moscow State University;
\\
e-mail: yurii.neretin@univie.ac.at
\\
URL:www.mat.univie.ac.at/$\sim$neretin

\end{document}